\documentclass[12pt,a4]{elsarticle}
\usepackage{amssymb}
\usepackage{pifont}
\usepackage{amsfonts}
\usepackage{amsmath}
\usepackage{graphicx}
\usepackage{setspace}
\usepackage[left=2cm,top=1.75cm,right=2cm]{geometry}

\newtheorem{lemma}{Lemma}
\newtheorem{remark}{Remark}
\usepackage[mathlines,pagewise]{lineno}
\newtheorem{theorem}{Theorem}

\numberwithin{equation}{section}

\begin{document}
\journal{Springer}
\begin{frontmatter}
\title{A short note on $q$-analogue of modified Stancu-Beta operators}

\author[]{Preeti Sharma Joshi}
\ead{preeti.iitan@gmail.com}

\author[]{Ghanshyam Singh Rathore}
\ead{icpam2017@gmail.com, ghanshyamsrathore@yahoo.co.in}

\address{Department of Mathematics \& Statistics,
Mohanlal Sukhadia University, Udaipur (Rajasthan), India}


\begin{abstract}
This paper deals with the modified $q$-Stancu-Beta operators and we have investigated the statistical approximation theorems for these operators with the help of the Korovkin type approximation theorem. We have also established the rates of statistical convergence by means of the modulus of continuity and the Lipschitz type maximal function. Our results show that the rates of convergence of our operators are at least as fast as the classical Stancu-Beta operators.
\end{abstract}
\end{frontmatter}

\section{Introduction}
In 1995, Stancu \cite{st321} introduced Beta operators $L_{n}$ of second
kind in order to approximate the Lebesgue integrable functions on the
interval $I=(0,\infty )$ as follows:
\begin{equation}\label{op1}
L_{n}(f;x)=\frac{1}{B(nx,n+1)}\int_{0}^{\infty }\frac{t^{nx-1}}{(1+t)^{nx+n+1}}%
f(t)dt,
\end{equation}
where, $x\in I, n\in \mathbb{N}$, $f$ is a real-valued functions defined on $I$, and $B$ is the Euler Beta function defined as:
\begin{equation*}
B(t,s)= 
\int_{0}^{\infty}\frac{x^{t-1}}{(1+x)^{t+s}}dt, \,\, t,s>0.
\end{equation*}
\newline
In 2012, Aral and Gupta \cite{aral2012q} introduced the $q$-analogue of Stancu-Beta operators as follows:%

\begin{equation}\label{op2}
L_{n}^{q}(f;x)=\frac{K(A,[n]_{q}x)}{B_{q}([n]_{q}x,[n]_{q}+1)}%
\int_{0}^{\infty /A}\frac{u^{[n]_{q}x-1}}{(1+u)_{q}^{[n]_{q}x+[n]_{q}+1}}%
f(q^{[n]_{q}x}u)d_{q}u,  \text{ } \forall n\in  \mathbb{N}
\end{equation}
where $q\in (0,\infty),$ and $x\in [0,\infty)$.
In this paper they estimated moments, and direct results in terms of modulus of continuity and their asymptotic formulae.
\newline
The following theorem on operators (\ref{op2}) is given by Aral and Gupta \cite{aral2012q}:
\begin{theorem}
Let $q=(q_{n})$ satisfy $0<q_{n}<1$ with $%
\lim\limits_{n\rightarrow \infty }q_{n}=1.$ For each $f\in C_{x^{2}}^{\ast
}[0,\infty )$, we have
\begin{equation*}
\lim_{n\rightarrow \infty }\Vert (L_{n}^{q_{n}}(f);.)-f\Vert _{x^{2}}=0,
\end{equation*}
where $C_{x^{2}}[0,\infty )$ denotes the subspace of all continuous
functions on $[0,\infty )$ such that $|f(x)|\leq M_{f},$ and $%
C_{x^{2}}^{\ast }[0,\infty )$ denotes the spaces of all $f\in
C_{x^{2}}[0,\infty )$ such that $\lim\limits_{x\rightarrow \infty }\frac{f(x)%
}{1+x^{2}}$ is finite. The norm on $C_{x^{2}}^{\ast }[0,\infty )$ is given
by $\Vert f\Vert _{x^{2}}=\sup\limits_{x\in [0,\infty )}\frac{\mid f(x)\mid
}{1+x^{2}}.$
\end{theorem}
After that Mursaleen and Khan \cite{mursaleen2013statistical} defined modified $q$-Stancu-Beta operators as follows:

\begin{equation}\label{op3}
L_{n}^{\ast }(f;q,x)=q~~\frac{K(A,[n]_{q}x)}{B_{q}([n]_{q}x,[n]_{q}+1)}%
\int_{0}^{\infty /A}\frac{u^{[n]_{q}x-1}}{(1+u)_{q}^{[n]_{q}x+[n]_{q}+1}}%
f(q^{[n]_{q}x}u)d_{q}u
\end{equation}%
\newline
where $x\geq 0,0<q\leq 1$. It is easy to verify that if $q=1$, then these
operators turns into the classical Stancu-Beta operators. \newline

In 2014, Cai \cite{cai2014approximation} introduced a new kind of modification of $q$-Stancu-Beta operators which preserve $x^2$ based on the concept of $q$-integer

\begin{equation}\label{op3}
\mathfrak{L}_{n,q}(f;x)= \frac{K(A,[n]_{q} v_{n}(x)}{B_{q}([n]_{q}v_{n}(x);[n]_{q}+1)}%
\int_{0}^{\infty /A}\frac{u^{[n]_{q}v_{n}(x)-1}}{(1+u)_{q}^{[n]_{q}v_{n}(x)+[n]_{q}+1}}%
f(q^{[n]_{q}v_{n}(x)}u)d_{q}u
\end{equation}%
where\\
$v_{n}(x)= \sqrt{ \frac{q[n]_{q} - q}{[n]_q} x^{2} +\frac{1}{4[n]_q^2}} -\frac{1}{2[n]_q},$ for $f \in C[0,\infty),~ n\in  \mathbb{N}, q\in (0,1)$   

\section{Some preliminary results}
In this section we give the following lemmas, which we need to prove our theorems.
\begin{lemma}\cite{aral2012q}
The following equalities hold:
\begin{equation*}
L_{n}^{q}(1;x)=1,~~L_{n}^{q}(t;x)= x \text{ and }~ L_{n}^{q}(t^{2};x)=\frac{([n]_{q}x+1)x}{q([n]_{q}-1)}.
\end{equation*}%
\end{lemma}

\begin{lemma}\cite{mursaleen2013statistical}
Note that $L_{n}^{\ast }(f;q,x)=L_{n}^{q}(f;x)\ $ and
from the Lemma 1 of Aral and Gupta \cite{aral2012q}, we have $%
L_{n}^{q}(1;x)=1,L_{n}^{q}(t;x)=x,L_{n}^{q}(t^{2};x)=\frac{([n]_{q}x+1)x}{%
q([n]_{q}-1)}.$ Hence for $x\geq 0,~0<q\leq 1$, we have%
\begin{equation*}
L_{n}^{\ast }(1;q,x)=q,~~L_{n}^{\ast }(t;q,x)=q \,x\text{ and }~L_{n}^{\ast
}(t^{2};q,x)=\frac{([n]_{q}x+1)x}{([n]_{q}-1)}.
\end{equation*}%
\end{lemma}

\begin{lemma}\cite{cai2014approximation}\label{cai1}
Let $q\in (0,1),~ x \in [0,\infty)$, we have\\
$(i) \mathfrak{L}_{n,q}(1;x) = 1\\
(ii) \mathfrak{L}_{n,q}(t;x) = \sqrt{ \frac{q[n]_{q} - q}{[n]_q} x^{2} +\frac{1}{4[n]_q^2}} -\frac{1}{2[n]_q}\\
(iii) \mathfrak{L}_{n,q}(t^2;x) = x^2$
\end{lemma}

\begin{remark}
Let $q\in (0,1)$, then for $x\in \lbrack 0,\infty )$, we
have
\begin{equation*}
\alpha_{n}(x) = \mathfrak{L}_{n,q}\left( t-x;x\right) =0
\end{equation*}%
and

\begin{equation*}
\delta_{n}(x) = \mathfrak{L}_{n,q}\left( {(}t-x{)}^{2};x\right) = 2x^2 - 2 x \sqrt{ \frac{q[n]_{q} - q}{[n]_q} x^{2} +\frac{1}{4[n]_q^2}} + \frac{x}{[n]_q}.
\end{equation*}%
\end{remark}
First, we recall certain notations of $q$-calculus and the details on $q$%
-integers can be found in \cite{er, kac2001quantum}.
\newline
The $q$ -integer $[k]_{q}$ for each nonnegative integer is defined as
\begin{equation*}
\lbrack k]_{q}:=\left\{
\begin{array}{ll}
\frac{(1-q^{k})}{(1-q)},~~~~~~q\neq 1 &  \\
~k,~~~~~~~~~~~q=1 &
\end{array}%
\right.
\end{equation*}%
\newline
\begin{equation*}
\lbrack k]_{q}!:=\left\{
\begin{array}{ll}
\lbrack k]_{q}[k-1]_{q}...[1]_{q},~~~~~~k\geq 1 &  \\
~1,~~~~~~~~~~~~~~~~~~~~~~~~~k=0 &
\end{array}%
\right.
\end{equation*}%
and
\begin{equation*}
\left[
\begin{array}{c}
n \\
k%
\end{array}%
\right] _{q}:=\frac{[n]_{q}!}{[k]_{q}![n-k]_{q}!}.
\end{equation*}%
The $q$-improper integral (Koornwinder \cite{koor312}) is defined as 
\begin{equation*}
\int_{0}^{\infty /A}f(x)d_{q}x=(1-q)\sum\limits_{n=-\infty }^{\infty }f\big(%
\frac{q^{n}}{A}\big)\frac{q^{n}}{A},~~~~A>0.
\end{equation*}%
%
\newline
Now, we extend the earlier work done by Cai \cite{cai2014approximation}. By considering the operators defined in \cite{cai2014approximation}, we shall obtain some approximation properties of Modified $q$ Stancu-Beta operators. We shall also estimate the rate of statistical convergence of these sequence of the operators. 

\section{Korovkin type statistical approximation properties}
%
Firstly, we recall the concept of statistical convergence for sequences of real numbers which was introduced by Fast \cite{fast1951convergence} and further studied by \cite{Fridy1985, Fridy1997, gadjiev2002some}  and many others.
\newline
Let $K\subseteq \mathbb{N}$ and $K_{n}=\left\{ j\leq n:j\in
K\right\}.$ Then the $natural~density$ of $K$ is defined by $\delta
(K)={\lim\limits_{n}}~ n^{-1}|K_{n}|$ if the limit exists, where $|K_{n}|$ denotes the cardinality of the set $K_{n}$.
\newline
 A sequence $x=(x_{j})_{j\geq1}$ of real numbers is said to be $%
statistically$ $convergent$ to $L$ provided that for every $\epsilon >0$ the
set $\{j\in \mathbb{N}:|x_{j}-L|\geq \epsilon \}$ has natural density zero,
i.e. for each $\epsilon >0$,
\begin{equation*}
\lim\limits_{n}\frac{1}{n}|\{j\leq n:|x_{j}-L|\geq \epsilon \}|=0.
\end{equation*}
It is denoted by $st-\lim\limits_{n}x_{n}=L$.
\newline
Doğru and Kanat \cite{dougru2012statistical}, defined the Kantorovich-type modification of Lupa\c{s} operators as follows:

\begin{equation}\label{dk}
\tilde{R}_n(f;q;x)=[n+1]\sum_{k=0}^{n}\bigg(\int_{\frac{[k]}{[n+1]}}^{\frac{[k+1]}{[n+1]}} ~ f(t) d_{q}t\bigg)\left(\begin{array}{c}n \\k \end{array}\right) \frac{q^{-k}q^{k(k-1)/2} x^k (1-x)^{(n-k)}}{(1-x+qx)\cdots(1-x+q^{n-1} x)}.
\end{equation}
They proved the following statistical Korovkin-type approximation theorem for operators (\ref{dk}).

\begin{theorem}
Let $q:=(q_n),~ 0 < q < 1$, be a sequence satisfying the following conditions:\\
\begin{equation}\label{thm1}
 st-\lim_n q_n = 1,~ st - \lim_n q_n^n = a ~( a < 1)~ and~ st - \lim_n \frac{1}{[n]_q} = 0,
\end{equation}
then if $f$ is any monotone increasing function defined on $[0, 1]$, for the positive linear operator $\tilde{R}_n(f;q;x)$, then
$$ st- \lim_n {\Vert \tilde{R}_n(f;q;\cdot) - f \Vert}_{C[0, 1]} = 0$$ holds. 
\end{theorem}
Do\u{g}ru \cite{do36} gave some examples so that $(q_{n})$ is statistically convergent to $1$ but it may not convergent to $1$ in the ordinary case.\\
\newline
Now, we consider a sequence $q=(q_{n}),$ $q_{n}\in $ $(0,1),$ such that

\begin{equation}\label{a1}
\lim\limits_{n\rightarrow \infty }q_{n}=1.
\end{equation}
The condition (\ref{a1}) guarantees that
$[n]_{q_{n}}\rightarrow
\infty $ as $n\rightarrow \infty .$

\begin{theorem} \label{thm2}
Let ${\mathfrak{L}_{n,q_n}}$ be the sequence of the operators (\ref{op3}) and the sequence
$q = (q_n)$ satisfies (\ref{thm1}). Then for any function $f \in C[0,\nu] \subset C[0,\infty),~ \nu > 0 ,$ we have
\begin{equation}\label{eq1thm2}
st- \lim_n \Vert {\mathfrak{L}_{n,q_n}}( f, \cdot) - f \Vert = 0,
\end{equation}
where $C[0,\nu]$ denotes the space of all real bounded functions $f$ which are continuous in $[0,\nu].$
\end{theorem}
\textbf{Proof.} ~~Let $f_{i}= t^i,$ where $i=0,1,2.$ 
As $\mathfrak{L}_{n,q_n}(1;x)=1,$ and $\mathfrak{L}_{n,q_n}(t^2;x)= x^2,$ (see Lemma (\ref{cai1}),~(\ref{eq1thm2}) holds true for i=0 and i=2.
So, $\mathfrak{L}_{n,q}(1;x)=1,$ it is clear that\\
$st-\lim\limits_{n}\|\mathfrak{L}_{n,q}(1;x)-1\|=0.$ and $st-\lim\limits_{n}\|\mathfrak{L}_{n,q}(t^2;x)- x^2\|=0.$
\newline
Finally, for i=1, by Lemma (\ref{cai1})(ii), we have

\begin{equation*}
\Vert \mathfrak{L}_{n,q_n}(t,x)-x\Vert
=\left\| \sqrt{ \frac{q[n]_{q} - q}{[n]_q} x^{2} +\frac{1}{4[n]_q^2}} -\frac{1}{2[n]_q}-x\right\| \leq  \bigg(1- \sqrt{ \frac{q[n]_{q} - q}{[n]_q}}\bigg) x +\frac{1}{2[n]_q}.
\end{equation*}
For given $\epsilon >0$, we define the following sets:
\begin{equation*}
U=\{k:\Vert{\mathfrak{L}}_{n,q_n}(t,x)-x\Vert\geq \epsilon\},
\end{equation*}
and 

\begin{equation}\label{e}
U'=\left\{k: \bigg(1- \sqrt{ \frac{q[k]_{q} - q}{[k]_q}}\bigg) x +\frac{1}{2[k]_q}\geq \epsilon\right\}.
\end{equation}
It is obvious that $U \subset U^{\prime }$, it can be written as

\begin{equation*}
\delta \left( \{k\leq n:\Vert \mathfrak{L}_{n,q_n}(t,x) - x \Vert \geq
\epsilon \}\right) \leq \delta \left( \{k \leq n: ~ \bigg(1- \sqrt{ \frac{q[k]_{q} - q}{[k]_q}}\bigg) x +\frac{1}{2[k]_q} \geq \epsilon \}\right).
\end{equation*}
By using (\ref{thm1}), we get

$$st-\lim_{n}\left(\bigg(1- \sqrt{ \frac{q[n]_{q} - q}{[n]_q}}\bigg) x +\frac{1}{2[n]_q}\right)=0.$$ 
So, we have

$$\delta \bigg(\big\{k \leq n: ~ \bigg(1- \sqrt{ \frac{q[n]_{q} - q}{[n]_q}}\bigg) x +\frac{1}{2[n]_q} \geq \epsilon \big\}\bigg) = 0,$$
then
$$st-\lim_{n}\|\mathfrak{L}_{n,q}(t,x)-x\|=0.$$
This completes the proof of theorem. 
\qed

\section{Weighted statistical approximation}
In this section, we obtain the Korovkin type weighted
statistical approximation by the operators defined in (\ref{op3}).
A real function $\rho $ is called a weight function if it is continuous on $%
\mathbb{R}$ and $\lim\limits_{\mid x\mid \rightarrow \infty }\rho (x)=\infty
,~\rho (x)\geq 1$ for all $x\in \mathbb{R}$.

Let by $B_{\rho }(\mathbb{R})$  denote the weighted space of real-valued
functions $f$ defined on $\mathbb{R}$ with the property $\mid f(x)\mid \leq
M_{f}~\rho (x)$ for all $x\in \mathbb{R}$, where $M_{f}$ is a constant
depending on the function $f$. We also consider the weighted subspace $%
C_{\rho }(\mathbb{R})$ of $B_{\rho }(\mathbb{R})$ given by $C_{\rho }(%
\mathbb{R})=\{f\in B_{\rho }(\mathbb{R}){:}$ $f$ continuous on $\mathbb{R}%
\} $. Note that $B_{\rho }(\mathbb{R})$ and $C_{\rho }(\mathbb{R})$ are
Banach spaces with $\Vert f\Vert _{\rho }=\sup\limits_{x\in R}\frac{\mid
f(x)\mid }{\rho (x)}.$ In case of weight function $~\rho _{0}=1+x^{2},$ we have $%
\Vert f\Vert _{\rho _{0}}=\sup\limits_{x\in R}\dfrac{\mid f(x)\mid }{1+x^{2}}%
.$

Now we are ready to prove our main result as follows:\newline
\begin{theorem}\label{thm3}
Let $\mathfrak{L}_{n,q_n}(f;x)$ be the sequence of the operators (\ref{op3})and the sequence $q=(q_{n})$ satisfies (\ref{thm1}). Then for any function $f\in C_{B}[0,\infty ),$ we have

\begin{equation*}
st-\lim_{n\rightarrow \infty }{\Vert }\mathfrak{L}_{n,q_n}(f;x) - f{\Vert }%
_{\rho _{0}}=0.
\end{equation*}
\end{theorem}
\textbf{Proof.} By Lemma (\ref{cai1})(iii), we have  $\mathfrak{L}_{n,q_n}(t^2, x) \leq Cx^2,$ where $C$ is a positive constant, $\mathfrak{L}_{n,q_n}(f;x)$ is a sequence of positive linear operator acting from $C_{\rho}[0,\infty)$ to $B_{\rho}[0,\infty)$. \\

Let $f_{i}= t^i,$ where $i=0,1,2.$ 
Since $\mathfrak{L}_{n,q_n}(1;x)=1,$ and $\mathfrak{L}_{n,q_n}(t^2;x)= x^2,$ (see Lemma (\ref{cai1}),~(\ref{eq1thm2}) holds true for i=0 and i=2.
So, $\mathfrak{L}_{n,q}(1;x)=1,$ it is clear that\\
$st-\lim\limits_{n}\|\mathfrak{L}_{n,q}(1;x)-1\|_{\rho_{0}}=0.$ and $st-\lim\limits_{n}\|\mathfrak{L}_{n,q}(t^2;x)- x^2\|_{\rho_{0}}=0.$
\newline
Finally, for i=1, by Lemma (\ref{cai1})(ii), we have

\begin{eqnarray*}
\| \mathfrak{L}_{n,q_n}(t,x)-x\|_{\rho_{0}} &=& \sup_{x\in [0,\infty)} \frac{|{\mathfrak{L}_{n,q_n}(t,x)-x|}}{1+x^2}\\ &\leq & \bigg(1- \sqrt{ \frac{q[n]_{q} - q}{[n]_q}}\bigg) \sup_{x\in [0,\infty)}\frac{x}{1+x^2}  +\frac{1}{2[n]_q} \sup_{x\in [0,\infty)}\frac{1}{1+x^2}\\
&\leq &\bigg(1- \sqrt{ \frac{q[n]_{q} - q}{[n]_q}}\bigg) +\frac{1}{2[n]_q}.
\end{eqnarray*}
Using (\ref{thm1}), we get

 $$st-\lim_{n}\left( \bigg(1- \sqrt{ \frac{q[n]_{q} - q}{[n]_q}}\bigg)  +\frac{1}{2[n]_q}\right)=0,$$
then $$st-\lim_{n}\|\mathfrak{L}_{n,q_n}(t,x)-x\|_{\rho_{0}}=0.$$
This completes the proof of the theorem.\qed


\section{Rates of statistical convergence}
In this section, by using the modulus of continuity, we will study rates of statistical convergence of operator (\ref{op3}) and Lipschitz type maximal functions are introduced.

\begin{lemma}\label{l3}
Let $0<q<1$ and $a\in[0,bq],~ b>0.$ The inequality
\begin{equation}
\int_{a}^{b}|t-x|d_qt\leq \left( \int_{a}^{b}|t-x|^2d_qt\right)^{1/2} \left(\int_{a}^{b} d_qt \right)^{1/2}
\end{equation}
is satisfied.
\end{lemma}

Let $C_B[0,\infty),$ the space of all bounded and continuous functions on $[0,\infty)$ and $x \geq 0.$
Then, for $\delta>0,$ the modulus of continuity of $f$ denoted by $\omega(f;\delta )$ is defined to be
\begin{equation*}
\omega(f;\delta)= \sup_{|{t- x}|\leq {\delta}}|f(t)-f(x)|,~t\in[0,\infty).
\end{equation*}
It is known that $\lim\limits_{\delta\rightarrow 0}\omega(f ; \delta) = 0$ for $f\in C_B[0,\infty)$ and also, for any $\delta > 0$ and each $t,x\geq 0,$ we have

\begin{equation}\label{2.2}
|f(t)-f(x)|\leq \omega(f;\delta)\left(1+\frac{|t-x|}{\delta}\right).
\end{equation}

\begin{theorem}
Let $(q_n)$ be a sequence satisfying (\ref{thm1}). For every non-decreasing $f\in C_B[0,\infty), ~x\geq 0$ and $n\in \mathbb{N},$ we have
$$|\mathfrak{L}_{n,q_n}(f,x)-f(x)| \leq 2 \omega(f;\sqrt{\delta_n(x)}),$$
where
\begin{eqnarray*}\label{d1}
\delta_n{(x)}= 2x^2 - 2 x \sqrt{ \frac{q[n]_{q} - q}{[n]_q} x^{2} +\frac{1}{4[n]_q^2}} + \frac{x}{[n]_q}.
\end{eqnarray*}
\end{theorem}
\textbf{Proof.}
Let  $f\in C_B[0,\infty)$ be a non-decreasing function and $x\geq 0$. Using linearity and positivity of the operators
 $\mathfrak{L}_{n,q_n}$ and then applying (\ref{2.2}), we get for $\delta > 0$ 
\begin{eqnarray*}
|\mathfrak{L}_{n,q_n}(f,x)-f(x)| &\leq & \mathfrak{L}_{n,q_n}\big(|f(t)-f(x)|,x\big)\\
&\leq & \omega(f,\delta)\big\{\mathfrak{L}_{n,q_n}(1,x)+ 
\frac{1}{\delta}{\mathfrak{L}_{n,q_n}}(|t-x|,x) \big\}.
\end{eqnarray*}
Taking
$\mathfrak{L}_{n,q_n}(1,x) = 1 $ and using Cauchy-Schwartz inequality,  we have
\begin{align*}
|\mathfrak{L}_{n,q_n}(f,x)-f(x)|
&\leq  \omega(f;\delta)\bigg\{1+\frac{1}{\delta}{\bigg(\mathfrak{L}_{n,q_n}((t-x)^2,x\big)}^{1/2}{\mathfrak{L}_{n,q_n}(1,x)}^{1/2}\bigg)\bigg\}\\
&\leq  \omega(f;\delta)\bigg[1+\frac{1}{\delta}\bigg\{2x^2 - 2 x \sqrt{ \frac{q[n]_{q} - q}{[n]_q} x^{2} +\frac{1}{4[n]_q^2}} + \frac{x}{[n]_q}\bigg\}^{1/2}\bigg].
\end{align*}
Taking $q = (q_n),$ a sequence satisfying (\ref{thm1}) and choosing $\delta = \delta_n(x)$ as in (\ref{d1}), the theorem is proved.\qed \\


Now we will give an estimate concerning the rate of approximation by means
of Lipschitz type maximal functions.\\
 In \cite{lb313}, Lenze introduced a Lipschitz type maximal function as%
\begin{equation*}
{f}_{\alpha }(x,y)=\sup\limits_{t>0,t\neq x}\frac{\mid f(t)-f(x)\mid }{{%
\mid }t-x{\mid }^{\alpha }}.
\end{equation*}%
\newline
In \cite{adb24}, the Lipschitz type maximal function space on $E\subset \lbrack
0,\infty )$ is defined as follows%
\begin{equation*}
\tilde{V}_{\alpha,E }=\{f=\sup (1+x)^{\alpha }~ {f}_{\alpha }(x,y)\leq M%
\frac{1}{(1+y)^{\alpha }};x\geq 0~{and}~y~\in E\},
\end{equation*}%
where $f$ is bounded and continuous function on $[0,\infty )$, $M$ is a
positive constant and $0<\alpha \leq 1$. \newline
Also, let $d(x,E)$ be the distance between $x$ and $E,$ that is,
$$d(x,E)= \inf \{|x-y|; y \in E\}.$$

\begin{theorem}\label{thm4}
If $\mathfrak{L}_{n,q_n}$ be defined by (\ref{op3}), then for all $%
f\in \tilde{V}_{\alpha ,E}$
\begin{equation}\label{e2.4.3}
\mid \mathfrak{L}_{n,q_n}(f,x)-f(x)\mid \leq M(\delta _{n}^{\frac{\alpha }{2}%
} ~ + ~d(x,E)),
\end{equation}%
where
\begin{equation}\label{delta2}
\delta_n{(x)}=2x^2 - 2 x \sqrt{ \frac{q[n]_{q} - q}{[n]_q} x^{2} +\frac{1}{4[n]_q^2}} + \frac{x}{[n]_q}
\end{equation}%
\end{theorem}
\textbf{Proof.} Let $x_0\in \bar{E}$, where $\bar{E}$ denote the closure of the set $E$. Then
we have\newline
\begin{equation*}
\mid f(t)-f(x)\mid \leq \mid f(t)-f(x_{0})\mid +\mid f(x_{0})-f(x)\mid .
\end{equation*}%
Since $\mathfrak{L}_{n,q_n} $ is a positive and linear operator, $f\in \tilde{V}%
_{\alpha ,E}$ and using the above inequality

\begin{equation*}
\mid {\mathfrak{L}_{n,q_n}}(f,x)-f(x)\mid ~ \leq {\mathfrak{L}_{n,q_n}}(\mid
f(t)-f(x_{0})\mid ;q_{n};x) + (\mid f(x_{0})-f(x)\mid)\mathfrak{L}_{n,q_n}(1,x)
\end{equation*}
\begin{equation}\label{e2.4.5}
\leq M\left( \mathfrak{L}_{n,q_n}({\mid }t-x_{0}{\mid }^{\alpha };q_{n};x)+{\mid }%
x-x_{0}{\mid }^{\alpha }{\mathfrak{L}_{n,q_n}}(1;q_{n};x)\right) .
\end{equation}%
Therefore, we have
\begin{equation*}
{\mathfrak{L}_{n,q_n}}\left( {\mid }t-x_{0}{\mid }^{\alpha };q_{n};x\right) \leq
{\mathfrak{L}_{n,q_n}}({\mid }t-x{\mid }^{\alpha };q_{n};x)+{\mid }x-x_{0}{\mid }%
^{\alpha }{\mathfrak{L}_{n,q_n}}(1,x).
\end{equation*}%
Now, we take $p=\frac{2}{\alpha }$ and $q=\frac{2}{(2-\alpha) }$ and by using the H\"{o}lder's inequality, one can write

\begin{equation*}
{\mathfrak{L}_{n,q_n}}\left( {(}t-x)^{\alpha };q_{n};x\right) \leq {\mathfrak{L}_{n,q_n}}\left( {(}t-x{)}^{2},x\right) ^{\alpha/2}({\mathfrak{L}_{n,q_n}}(1,x){)}%
^{{(2-\alpha)}/{2}}
\end{equation*}%
\begin{equation*}
+{\mid }x-x_{0}{\mid }^{\alpha }{\mathfrak{L}_{n,q_n}}(1,x)
\end{equation*}
\begin{equation*}
=\delta _{n}^{\frac{\alpha }{2}} + {\mid }x-x_{0}{%
\mid }^{\alpha }.
\end{equation*}%
Substituting this in ($\ref{e2.4.5}$), we get ($\ref{e2.4.3}$).\newline
This completes the proof of the theorem.\qed

\section{Concluding remarks}
\begin{remark}
Observe that by the conditions in (\ref{thm1}),
$$st-\lim\limits_n \delta_n = 0.$$
By (\ref{delta2}), we have
$$st-\lim\limits_n\omega(f;\delta_n) = 0.$$
This gives us the pointwise rate of statistical convergence of the operators  $\mathfrak{L}_{n,q_n}(f,x)$ to $f(x).$
\end{remark}

\begin{remark}
If we take $E=[0,\infty )$ in Theorem \ref{thm4}, since $d(x,E)=0,$
then we get for every $f\in \tilde{V}_{\alpha ,[0,\infty )}$
\begin{equation*}
{\mid }\mathfrak{L}_{n,q_n}(f,x)-f(x){\mid }\leq M\delta _{n}^{\frac{%
\alpha }{2}}
\end{equation*}%
where $\delta _{n}$ is defined as in (\ref{delta2}).
\end{remark}


\hspace{-0.9cm}

\end{document}